\documentclass[11pt,a4paper]{article}
\usepackage{french,amsmath,isolatin1,amssymb,euscript,amscd,epsfig}

\title{Un théorème de Cartier-Milnor-Moore-Quillen pour les
 bigèbres dendriformes et les algèbres braces} 

\author{Frédéric Chapoton}

\date{}

\newcommand{\PP}{\mathcal{P}}
\newcommand{\QQ}{\mathcal{Q}}
\newcommand{\JJ}{\mathcal{J}}
\newcommand{\CC}{\mathcal{C}}
\newcommand{\UU}{\mathcal{U}}
\newcommand{\DD}{\mathcal{D}_+}

\newcommand{\ZZ}{\mathbb{Z}}
\newcommand{\RR}{\mathbb{R}}
\newcommand{\kk}{\mathbb{K}}
\newcommand{\smod}{\mathbb{S}}

\newcommand{\sym}{\mathfrak{S}}

\newcommand{\alg}{\operatorname{Alg}}
\newcommand{\vect}{\operatorname{Vect}_\kk}
\newcommand{\mont}{\nearrow}
\newcommand{\desc}{\searrow}
\newcommand{\plig}{\operatorname{Pli}}
\newcommand{\id}{\operatorname{Id}}
\newcommand{\prim}{\operatorname{Prim}}
\newcommand{\entr}{\operatorname{Entr}}
\newcommand{\aut}{\operatorname{End}}
\newcommand{\appe}{\mathcal{APE}}
\newcommand{\ang}{\operatorname{Angles}}
\newcommand{\cor}{\operatorname{Corl}}

\newcommand{\com}{\operatorname{Com}}
\newcommand{\prelie}{\operatorname{PreLie}}

\newcommand{\lie}{\operatorname{Lie}}
\newcommand{\ass}{\operatorname{As}}
\newcommand{\brac}{\operatorname{Brace}}

\newcommand{\dend}{\operatorname{Dend}}

\newcommand{\zin}{\operatorname{Zin}}
\newcommand{\ape}{\operatorname{APE}}

\newcommand{\viv}{\stackrel{v}{\vee}}
\renewcommand{\hom}{\operatorname{Hom}}

\newtheorem{theo}{Théorème}
\newtheorem{defi}{Définition}
\newtheorem{prop}{Proposition}
\newtheorem{coro}{Corollaire}
\newtheorem{lemm}{Lemme}
\newtheorem{sublemm}{Sous-lemme}

\newenvironment{preuve}{{\sc Preuve : ~}}{\hfill\ensuremath{\square}}

\begin{document}
\maketitle

\begin{abstract}
  Le théorème classique de Cartier-Milnor-Moore-Quillen donne une
  équivalence entre la catégorie des bigèbres cocommutatives connexes
  et la catégorie des algèbres de Lie. On définit ici une équivalence
  analogue entre la catégorie des bigèbres dendriformes connexes et la
  catégorie des algèbres braces. Cette équivalence est donnée par le
  foncteur ``Prim'' des éléments primitifs d'un bigèbre dendriforme et
  le foncteur ``algèbre dendriforme enveloppante'' d'une algèbre
  brace.
\end{abstract}



Classification AMS 2000 : 18D50, 16W30, 17B35

Mots-Clés : algèbre dendriforme, algèbre brace, cogèbre


\section*{Introduction}

On étudie ici un analogue du couple formé par une algèbre
de Lie et son algèbre associative enveloppante. Le rôle de la
notion d'algèbre associative est joué par celle d'algèbre
dendriforme, introduite par Loday
\cite{loday-cras,loday-overview,dialgebra} comme structure duale au
sens des opérades de celle de digèbre. Le rôle de la
notion d'algèbre de Lie est tenu par les algèbres braces,
apparues d'abord dans les articles de Kadeishvili \cite{kadeishvili88}
et Getzler \cite{getzler93} sur l'homologie des algèbres
$A_{\infty}$, et plus récemment dans les travaux sur le complexe
de Hochschild et la conjecture de Deligne, notamment ceux de
Kontsevich et Soibelman \cite{kont-soib} et de Gerstenhaber et Voronov
\cite{gerst-voro-higher,gerst-voro-imrn}. La notion d'algèbre
pré-Lie issue des travaux de Gerstenhaber \cite{gerst63} et
Vinberg \cite{vinberg63} joue aussi un rôle dans la suite.

On réinterprète en termes d'opérades des résultats de Marìa Ronco qui
montrent qu'une algèbre dendriforme est en particulier une algèbre
brace \cite{ronco}, en introduisant l'opérade des algèbres braces et
un morphisme d'opérades dans celle des algèbres dendriformes.

On introduit une notion de bigèbre dendriforme, voisine de la notion
d'algèbre de Hopf dendriforme introduite indépendamment par Ronco dans
\cite{ronco}, mais qui en diffère par la présence d'une unité. La
présence de cette unité offre des avantages techniques, dont celui de
simplifier certaines démonstrations.

On définit ensuite l'algèbre dendriforme enveloppante d'une algèbre
brace $B$, notée $\widetilde{\UU}(B)$. Par analogie avec le fait que
l'algèbre associative enveloppante d'une algèbre de Lie est une
bigèbre, on montre que l'algèbre dendriforme enveloppante d'une
algèbre brace est une bigèbre dendriforme.

Si $W$ est une bigèbre dendriforme, il est démontré dans \cite{ronco}
que l'ensemble $\prim(W)$ des éléments primitifs pour la structure de
cogèbre sous-jacente est une sous-algèbre brace de $W$.

On montre que, si $B$ est une algèbre brace, alors
$B=\prim\widetilde{\UU}(B)$ et que, réciproquement, si $D$ est une
bigèbre dendriforme connexe (comme cogèbre), alors $D\simeq
\widetilde{\UU}(\prim(D))$. On obtient ainsi une équivalence de
catégories qui forme un analogue dendriforme/brace du
théorème de Cartier-Milnor-Moore-Quillen.

\vspace{0.3cm}

Le plan de l'article est le suivant :

La première partie contient des rappels sur les différents
types d'algèbres utilisés et sur certaines relations entre
eux.

Dans la seconde partie, on introduit l'opérade $\brac$ qui
décrit les algèbres braces en termes d'arbres enracinés
plans, on en donne une présentation par générateurs et
relations et on définit un morphisme d'opérades de $\brac$
dans $\dend$. On retrouve ainsi le fait que toute algèbre
dendriforme est une algèbre brace, établi par \cite{ronco}. On
démontre enfin deux propositions techniques qui seront
utilisées dans la suite de l'article.

La troisième partie est consacrée à l'introduction des notions
d'algèbre dendriforme unitaire et de bigèbre dendriforme. On rappelle
que, d'après Ronco, les éléments primitifs d'une bigèbre dendriforme
forment une sous-algèbre brace, puis on énonce deux lemmes sur les
bigèbres dendriformes, qui joueront un rôle important dans la suite.
On montre ensuite que l'algèbre dendriforme libre est une bigèbre
dendriforme, pour le coproduit introduit par Loday et Ronco dans
\cite{loday-ronco}.

La quatrième partie contient des énoncés généraux sur les suites
exactes courtes d'opérades, les idéaux d'opérades et les adjoints des
foncteurs d'oubli entre catégories d'algèbres induits par les
morphismes d'opérades, pour lesquels on n'a pas trouvé de références
dans la littérature.

Dans la cinquième partie, on définit l'algèbre dendriforme
enveloppante d'une algèbre brace $L$, et l'on montre que
c'est une bigèbre dendriforme, au sens de la seconde partie. On
applique cette construction aux cas des algèbres braces libres et des
algèbres braces de produits nuls.

Enfin, dans la sixième partie, on montre que le foncteur
\textit{algèbre dendriforme enveloppante} $\widetilde{\UU}$
établit une équivalence de catégories entre la
catégorie des algèbres braces et celle des bigèbres
dendriformes connexes, dont le quasi-inverse (inverse au sens des
équivalences de catégories) est le foncteur $\prim$ des
éléments primitifs.

\section{Rappels sur les algèbres dendriformes et pré-Lie} 

On rappelle ici, pour la commodité du lecteur et pour les notations,
la définition des différents types d'algèbres utilisés. Les notions
d'algèbre dendriforme et d'algèbre de Leibniz duale ont été
introduites et développées par Loday, voir
\cite{loday-cras,loday-overview,dialgebra}. Les algèbres pré-Lie ont
été introduites indépendamment par Gerstenhaber \cite{gerst63} et
Vinberg \cite{vinberg63}. On renvoie aussi à \cite {rooted} pour des
résultats sur l'opérade correspondante.

\vspace{0.5cm}

Une algèbre dendriforme est un espace vectoriel $W$ muni de deux
produits notés $\prec$ et $\succ$ de $W\otimes W$ dans $W$
tels que
\begin{align}
  (x\prec y)\prec z&=x\prec (y\prec z)+x\prec (y\succ z),\\
  (x\succ y)\prec z&=x\succ (y\prec z),\label{pas-parenthese}\\
  x\succ (y\succ z)&=(x\succ y)\succ z+(x\prec y)\succ z.
\end{align}

On note $D(V)$ l'algèbre dendriforme libre sur un espace vectoriel
$V$, voir \cite[2.4]{dialgebra} pour sa définition en termes de
l'opérade $\dend$. Si $(W,\prec,\succ)$ est une algèbre
dendriforme, on note $*$ le produit défini par $x*y=x\prec
y+x\succ y$. C'est un produit associatif. On note $W_{\ass}$
l'algèbre associative $(W,*)$. Cette construction définit un
morphisme d'opérades $\ass\rightarrow \dend$.

\vspace{0.5cm}

Une algèbre pré-Lie est un espace vectoriel $W$ muni d'une opération
$\circ$ de $W\otimes W$ dans $W$ telle que
\begin{equation}
  (x\circ y)\circ z- x\circ (y\circ z)=(x\circ z)\circ y-x\circ (z\circ y).
\end{equation}
  
On note $PL(V)$ l'algèbre pré-Lie libre sur un espace vectoriel $V$,
voir \cite{rooted} pour sa définition en termes de l'opérade
$\prelie$. Si $(W,\circ)$ est une algèbre pré-Lie, on note $[\, ,\,]$ le
produit défini par $[x,y]=x\circ y-y\circ x$. C'est un crochet de Lie.
On note $W_{\lie}$ l'algèbre de Lie $(W,[\, ,\,])$. Cette construction
définit un morphisme d'opérades $\lie\rightarrow \prelie$.

\vspace{0.5cm}

Une algèbre de Leibniz duale ou $\zin$-algèbre est un espace
vectoriel $W$ muni d'un produit noté $\cdot$ de $W\otimes W$ dans
$W$ tel que
\begin{equation}
  (x\cdot y)\cdot z=x\cdot (y\cdot z)+x\cdot (z\cdot y).
\end{equation}

On note $Z(V)$ la $\zin$-algèbre libre sur un espace vectoriel
$V$, voir \cite[7.1]{dialgebra} pour sa définition en termes de
l'opérade $\zin$. Si $(W,\,\cdot\,)$ est une $\zin$-algèbre, on
note $*$ le produit défini par $x*y=x\cdot y+y\cdot x$. C'est un
produit associatif et commutatif. On note $W_{\com}$ l'algèbre
associative commutative $(W,*)$. Cette construction définit un
morphisme d'opérades $\com\rightarrow \zin$.

\vspace{0.5cm}

Si $(W,\prec,\succ)$ est une algèbre dendriforme, on note $\circ$
le produit défini par $x\circ y=x\prec y-y\succ x$. C'est un
produit pré-Lie. On note $W_{\prelie}$ l'algèbre pré-Lie
$(W,\circ)$. Cette construction définit un morphisme d'opérades
$\prelie\rightarrow \dend$.

\vspace{0.5cm}

En posant $y\succ x=x\prec y=x\cdot y$, on peut considérer une
$\zin$-algèbre $W$ comme une algèbre dendriforme, qui est
symétrique au sens où $x\prec y=y\succ x$, pour tous $x,y$ dans
$W$. Ceci définit un morphisme d'opérades $\dend\rightarrow
\zin$.

Réciproquement, on observe que toute algèbre dendriforme
symétrique au sens ci-dessus provient d'une $\zin$-algèbre par
cette construction. Il y a donc équivalence entre la notion
d'algèbre dendriforme symétrique et celle de $\zin$-algèbre. La
condition de symétrie entraîne la commutativité de l'algèbre
associative $W_{\ass}$.

\vspace{0.5cm}

On rappelle les morphismes d'opérades bien connus $\lie\to
\ass\to\com$. Les foncteurs correspondants entre catégories d'algèbres
consistent respectivement à associer à une algèbre associative
l'algèbre de Lie dont le crochet est le commutateur du produit
associatif et à considérer une algèbre commutative comme une algèbre
associative.

En résumé, on a le diagramme suivant de morphismes d'opérades :
\begin{equation}
\label{diagre}
\begin{CD}
\zin @((( \dend @((( \prelie\\
@AAA @AAA @AAA \\
\com @((( \ass  @((( \lie.\\
\end{CD}
\end{equation}
On vérifie aisément qu'il est commutatif.

\vspace{0.3cm}

Si $\PP=(\PP(k))_{k\geq 1}$ est une opérade telle que $\PP(1)$ soit
engendré par l'identité, alors il existe un morphisme naturel de
$\PP$ dans l'opérade triviale (réduite à l'identité). On note
$\PP_+=(\PP(k))_{k\geq 2}$ le noyau de ce morphisme d'augmentation.

Les deux lignes du diagramme (\ref{diagre}) sont exactes au sens
suivant : on a $\com\simeq \ass/<\lie_+>$ et
$\zin\simeq\dend/<\prelie_+>$, où $<\,>$ désigne l'idéal engendré.
Ceci signifie que les algèbres commutatives sont exactement les
algèbres associatives dont le crochet de Lie sous-jacent est nul. La
proposition similaire pour les $\zin$-algèbres est l'identification
décrite plus haut avec les algèbres dendriformes symétriques.

\vspace{0.3cm}

Si $\PP$ est une opérade, on notera $\PP-\alg$ la catégorie
des $\PP$-algèbres.

\section{L'opérade des algèbres braces}

On définit une opérade $\ape$ sur les arbres enracinés
plans, qui a été implicitement introduite par Kontsevich et
Soibelman \cite{oral-x,kont-soib}. On en donne ensuite une
présentation par générateurs et relations qui permet de
l'identifier avec l'opérade $\brac$ des algèbres braces
utilisées par Gerstenhaber et Voronov
\cite{gerst-voro-imrn,gerst-voro-higher}. Cette description en termes
d'arbres présente l'avantage de rendre visuelle la combinatoire des
produits braces. Les produits braces ont été généralisés
par Akman \cite{akman-multi,akman-master} en des produits
multi-braces, dont l'éventuelle relation avec ce qui suit reste à
clarifier. On montre que le morphisme $\prelie\to\dend$ se factorise
via l'opérade $\brac$, le morphisme $\brac\to\dend$ étant une
interprétation en termes d'opérades des résultats de
\cite{ronco}. La section se termine par des propositions sur les
idéaux de $\dend$ engendrés par les images de $\brac_+$ et
$\prelie_+$.

\subsection{Définition}
Soit $I$ un ensemble fini non vide. On appelle arbre enraciné plan sur
$I$ la donnée d'un graphe connexe sans boucles sur l'ensemble de
sommets $I$, muni d'un sommet distingué appelé la racine et d'un
plongement dans le plan. On dessine les arbres avec la racine en bas
et on considère les arêtes comme orientées vers la racine.

\begin{figure}[htbp]
  \begin{center}
    \leavevmode
    \epsfig{file=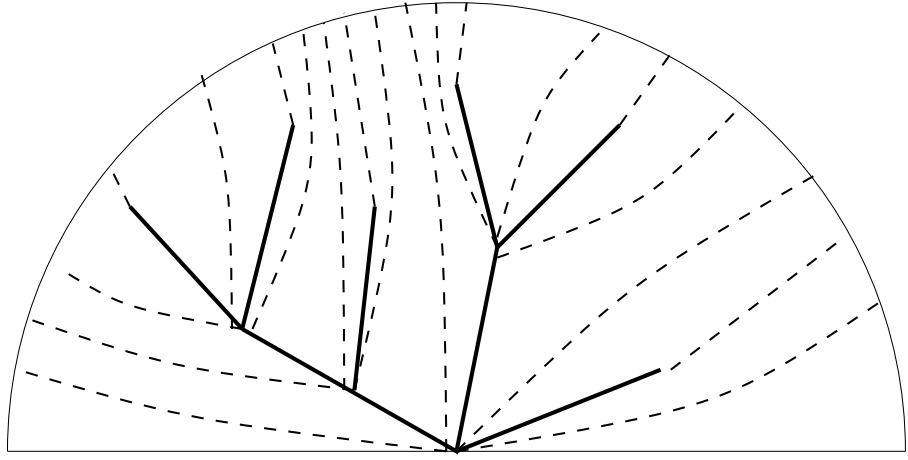,width=12cm}
    \caption{Angles d'un arbre plan}\label{figangle}
  \end{center}
\end{figure}

Supposons l'arbre $T$ dessiné dans le demi-disque supérieur ouvert
\begin{equation*}
  \DD=\{(x,y)\in\RR^2\mid y>0\text{ et }x^2+y^2<1\},
\end{equation*}
sauf la racine placée au point $x=y=0$. On appelle angle de $T$
une paire $(s,\alpha)$ où $s$ est un sommet de $T$ et $\alpha$ une
composante connexe de $B_{\epsilon}(s)\cap(\DD\setminus T)$ où
$B_{\epsilon}(s)$ est un petit disque de centre $s$. On note $\ang(T)$
l'ensemble des angles de $T$. On munit naturellement $\ang(T)$ d'un
ordre total de gauche à droite, de la façon suivante. En
considérant un angle comme une direction issue d'un sommet, on
peut tracer un chemin de chaque angle vers un point du cercle
unité. On peut aussi supposer que ces chemins ne se coupent pas.
On obtient alors un point du demi-cercle associé à chaque
angle. L'ordre de ces points de gauche à droite ne dépend pas
des chemins choisis pourvu qu'ils ne se coupent pas. Ceci définit
l'ordre total voulu, voir la figure \ref{figangle}.

Cette notion d'angle est due à Kontsevich \cite{oral-x}, \cite[p.
29-30]{kont-soib}, de même que la composition de l'opérade qui suit.

On note $\appe(I)$ l'ensemble des arbres enracinés plans sur $I$ et
$\ape(I)$ le $\ZZ$-module libre sur $\appe(I)$. Soit $T\in\appe(I)$. Si
$j$ est un sommet de $T$, on note $\entr(T,j)$ l'ensemble totalement
ordonné de gauche à droite des arêtes entrantes en $j$. On peut alors
considérer l'ensemble des fonctions croissantes au sens large de
$\entr(T,j)$ dans $\ang(S)$.

\begin{figure}[htbp]
  \begin{center}
   \epsfig{file=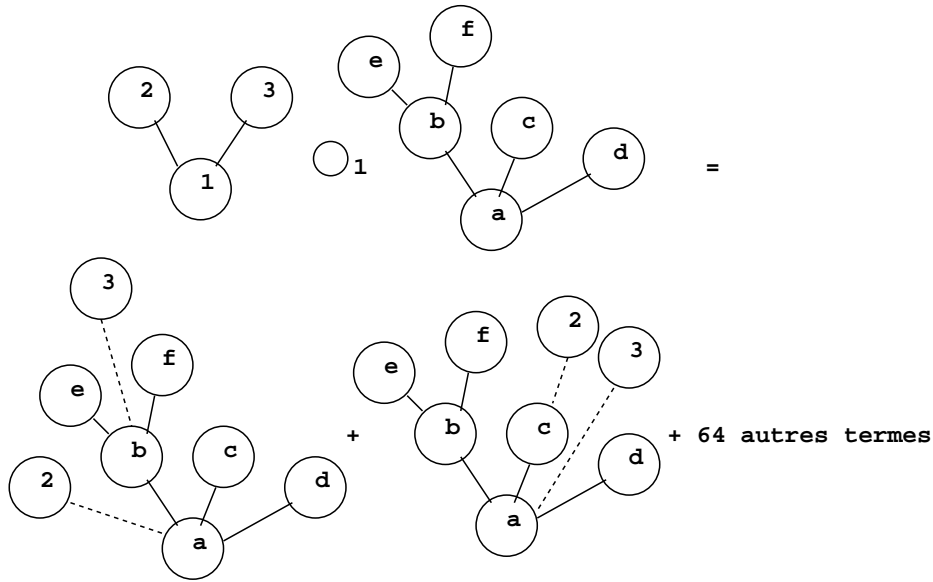,width=13cm}
    \caption{Quelques termes d'une composition}\label{fig-compo}
  \end{center}
\end{figure}

On définit une opérade $\ape$ sur les modules $\ape(I)$. La
composition d'un arbre $S\in\appe(I)$ au sommet $j$ d'un arbre $T\in
\appe(J)$ est définie par
\begin{equation}
  T\circ_j S=\sum_{f : \entr(T,j)\to \ang(S)} T\circ^f_j S,
\end{equation}
où $f$ est une application croissante au sens large et $ T\circ^f_j S$
est l'arbre obtenu par substitution de l'arbre $S$ au sommet $j$ de
$T$, les arêtes de $T$ entrantes dans $j$ étant greffées sur $S$ selon
l'application $f$. On donne un exemple dans la figure \ref{fig-compo},
avec $11$ angles pour $S$ et $2$ arêtes entrantes dans le sommet $1$
de $T$, ce qui donne $66$ termes dans la composition.

\vspace{0.5cm}

\begin{prop}
  La composition ci-dessus définit une opérade $\ape$.
\end{prop}
\begin{preuve}
  On vérifie aisément les axiomes d'unité, d'associativité et
  d'équivariance.
\end{preuve}

\subsection{Présentation}

On montre maintenant que l'opérade $\ape$ admet une présentation par
générateurs et relations qui permet de l'identifier à l'opérade des
algèbres braces.

\vspace{0.3cm}

On note $\cor_n(1;2,\dots,n+1)$ ou simplement $\cor_n$ la corolle à
$n$ feuilles indexée par $\{1,\dots,n+1\}$, la racine portant l'indice
$1$, les feuilles les indices $2,\dots,n+1$ dans l'ordre de gauche à
droite. On va prendre des éléments $C_n$ correspondant à ces corolles
comme générateurs. Les relations entre ces corolles sont celles entre
les opérations $x_1\{ x_2,\dots,x_{n+1}\}$ des algèbres braces, voir
\cite[Eq. (6)]{gerst-voro-higher} par exemple. Ces relations
expriment une composition de corolles à la racine comme une somme de
compositions multiples aux feuilles.

\begin{prop}
  L'opérade $\ape$ est isomorphe au quotient de l'opérade libre sur la
  famille de générateurs $C_n$ en $n+1$ variables pour $n\geq 1$ (on
  note $\{x_1 \mid x_2,\dots,x_{n+1}\}$ les opérations correspondantes)
  par les relations qui s'expriment en termes d'opérations de la façon
  suivante :
  \begin{multline}
    \label{relations}
    \{\{z\mid x_1,\dots,x_{n}\}\mid y_1,\dots,y_{m}\}=\\
    \sum \{z \mid Y_0,\{x_1\mid Y_1 \},Y_2,\{x_2\mid Y_3 \},Y_4
    \dots,Y_{2n-2},\{x_n\mid Y_{2n-1}\},Y_{2n}\},
  \end{multline}
  où la somme porte sur les partitions de l'ensemble ordonné
  $\{y_1,\dots,y_m\}$ en intervalles successifs éventuellement vides
  $Y_0\sqcup Y_1\sqcup\dots \sqcup Y_{2n}$.
\end{prop}
\begin{preuve}
  On note $\brac$ l'opérade quotient considérée. Elle décrit
  exactement les algèbres braces au sens de \cite{gerst-voro-higher}.
  On définit un morphisme d'opérades $\kappa$ de l'opérade libre sur
  les $C_n$ dans $\ape$ en posant $\kappa(C_n)=\cor_n$. Comme les
  corolles vérifient les relations (\ref{relations}), le morphisme
  $\kappa$ passe au quotient et définit un morphisme d'opérades $\beta
  : \brac\to\ape$.
  
  La surjectivité de $\kappa$ est claire, car tout arbre s'obtient par
  composition de corolles aux feuilles. Le morphisme $\beta$ est donc
  surjectif.
  
  D'autre part, on déduit des relations (\ref{relations}) que l'on peut
  définir un morphisme surjectif de $\ZZ$-modules $\gamma$ de
  $\ape(n)$ sur $\brac(n)$. La composée $\beta\circ \gamma$ est
  surjective, donc un isomorphisme de $\ZZ$-modules libres. Par
  conséquent $\gamma$ est injective, donc bijective et $\beta$ est un
  isomorphisme d'opérades.
\end{preuve}

Cette proposition permet d'identifier les opérades $\ape$ et $\brac$.
On appelle algèbres braces les $\brac$-algèbres. On note $B(V)$
l'algèbre brace libre sur un espace vectoriel $V$.

\subsection{Un morphisme de $\brac$ dans $\dend$}

On définit dans cette section un morphisme d'opérades de $\brac$ dans
$\dend$.

\vspace{0.3cm}

On définit par récurrence des éléments $\mont_n $ et $\desc^n $ de
$\dend$ :
\begin{equation}
  \mont_1(x_1)=x_1, 
  \mont_n(x_1,\dots,x_{n})=x_1\prec\,\mont_{n-1}(x_2,\dots,x_{n}),
\end{equation}
et 
\begin{equation}
  \desc^1(x_1)=x_1,
  \desc^{n+1}(x_1,\dots,x_{n+1})=\desc^n (x_1,\dots,x_n)\succ x_{n+1}.
\end{equation}
Par la suite, on notera simplement $\mont$ et $\desc$ en omettant le
nombre d'arguments lorsque cela ne prête pas à confusion.

\vspace{0.5cm}

On a alors la proposition suivante :

\begin{prop}
  Il existe un unique morphisme d'opérades $\psi$ de $\brac$ dans
  $\dend$ tel que pour tout $n$, $\psi(\cor_n)\in\dend(n+1)$ soit
  l'opération qui à $x_1,\dots,x_{n+1}$ associe
  \begin{equation}
    \sum_{i=1}^{n+1} 
    (-1)^i \mont(x_2,\dots,x_i) \succ x_1 \prec  
    \, \desc(x_{i+1},\dots,x_{n+1}),
  \end{equation}
  expression sans ambiguïté par la formule
  $(\ref{pas-parenthese})$, et où on convient que les termes
  pour $i=1$ et $i=n+1$ sont respectivement
  \begin{equation}
    x_1 \prec  
    \, \desc(x_{2},\dots,x_{n+1})\text{  et  }
    \mont(x_2,\dots,x_{n+1}) \succ x_1.
  \end{equation}
\end{prop}
\begin{preuve}
  L'unicité est claire puisque que les opérations $\cor_n$
  engendrent $\brac$.
  
  On peut prouver la compatibilité aux relations (\ref{relations}) en
  adaptant la preuve du théorème 3.4 de \cite{ronco}. On peut aussi en
  donner une preuve légèrement différente, qu'on ne détaillera pas
  ici, basée sur la description de l'opérade $\dend$ en termes
  d'arbres binaires plans \cite{dialgebra}.
\end{preuve}

\subsection{Un morphisme de $\prelie$ dans $\brac$}

On définit dans cette section un morphisme de $\prelie$ dans $\brac$.
On renvoie le lecteur à \cite{rooted} pour la description de l'opérade
$\prelie$ en termes d'arbres non-plans enracinés et la définition
précise de cette notion.

\begin{prop}
  Il existe un unique morphisme d'opérades $\phi$ de $\prelie$ dans
  $\brac$ tel que 
  \begin{equation}
    \phi(x_1\circ x_2)=\{x_1 \mid x_2\}.
  \end{equation}
  
  Ce morphisme est donné par la symétrisation des arbres : si $T$ est
  un arbre enraciné non-plan indexé par un ensemble $I$, $\phi(T)$ est
  la somme des arbres enracinés plans isomorphes à $T$ comme arbres
  enracinés non-plans. 
\end{prop}
\begin{preuve}
  L'unicité résulte du fait que $\prelie$ est engendrée par
  l'opération $x_1\circ x_2$. L'existence est un calcul immédiat avec
  la relation qui définit $\prelie$. Pour la seconde assertion, il
  suffit alors de vérifier que l'application de symétrisation
  définit bien un morphisme d'opérades, ce qui est clair par la forme
  des compositions en termes d'arbres dans les opérades $\prelie$ et
  $\brac$, voir section 2.1 et \cite{rooted}.
\end{preuve}

\vspace{0.3cm}

Il est clair que le morphisme composé $\psi\circ \phi$ coïncide avec
le morphisme de $\prelie\to\dend$ défini dans la première section.

\subsection{Idéal engendré par l'image de $\brac_+$ dans $\dend$}

La proposition suivante permet d'identifier l'opérade quotient de
$\dend$ par l'idéal engendré par $\psi(\brac_+)$.

\begin{prop}
  \label{quotient-zin}
  Les images de $\prelie_+$ et de $\brac_+$ dans $\dend$ engendrent le
  même idéal. L'opérade quotient de $\dend$ par cet idéal est
  l'opérade $\zin$.
\end{prop}
\begin{preuve}
  On a déjà vu dans la section 1 que $\zin$ est le quotient de $\dend$
  par l'idéal engendré par $\prelie_+$. Il reste donc à montrer la
  première assertion.
  
  Comme $\phi$ fait de $\prelie$ une sous-opérade de $\brac$, on a
  clairement une inclusion $\big{<}\psi\circ\phi(\prelie_+)
  \big{>}\subset \big{<}\psi(\brac_+) \big{>}$. Pour montrer la
  réciproque, il suffit de montrer que $\psi(\cor_n)$ appartient à
  $\big{<}\psi\circ\phi(\prelie) \big{>} $, pour tout $n\geq 1$.
  
  On rappelle que $\psi(\cor_n)$ est, par définition, l'opération qui
  à $x_1,\dots,x_{n+1}$ associe
  \begin{multline}
    - x_1 \prec \, \desc(x_{2},\dots,x_{n+1})\succ x_1
    +(-1)^{n+1}\mont(x_2,\dots,x_{n+1})
    \\+\sum_{i=2}^{n} (-1)^i
    \mont(x_2,\dots,x_i) \succ x_1 \prec \,
    \desc(x_{i+1},\dots,x_{n+1}).
  \end{multline}

  Pour chaque terme de la somme de $i=2$ à $i=n$, on utilise la relation
  \begin{multline*}
    \mont(x_2,\dots,x_i) 
    \succ x_1 \prec \,\desc(x_{i+1},\dots,x_{n+1})\equiv \\
    \big{(}\mont(x_2,\dots,x_i) \,* 
    \,\desc(x_{i+1},\dots,x_{n+1})\big{)}\succ x_1,
  \end{multline*}
  modulo $\big{<}\psi\circ\phi(\prelie) \big{>}$. On réécrit aussi le
  terme correspondant à $i=1$ :
  \begin{equation*}
    x_1 \prec \,\desc(x_{2},\dots,x_{n+1})\equiv 
    \desc(x_{2},\dots,x_{n+1})\succ x_1,
  \end{equation*}
  modulo $\big{<}\psi\circ\phi(\prelie) \big{>}$. On obtient donc au
  total
  \begin{multline*}
    -\desc(x_{2},\dots,x_{n+1})\succ x_1
    +(-1)^{n+1}\mont(x_2,\dots,x_{n+1}) \succ x_1\\
    +\sum_{i=2}^{n} (-1)^i \big{(}\mont(x_2,\dots,x_i) \,*
    \,\desc(x_{i+1},\dots,x_{n+1})\big{)}\succ x_1,
  \end{multline*}
  c'est-à-dire une expression de la forme $Z\succ x_1$, où
  \begin{multline*}
    Z= -\desc(x_{2},\dots,x_{n+1})+(-1)^{n+1}\mont(x_2,\dots,x_{n+1})\\
    +\sum_{i=2}^{n} (-1)^i \big{(}\mont(x_2,\dots,x_i) \,* 
    \,\desc(x_{i+1},\dots,x_{n+1})\big{)}.
  \end{multline*}
  Le lemme 2.6 de \cite{ronco} dit exactement que $Z$ est nul dans
  $\dend$, ce qui termine la démonstration.
\end{preuve}

\subsection{Idéal à gauche engendré par l'image de $\brac$}

On montre dans cette section que l'idéal à gauche (au sens défini
ci-dessous) engendré par l'image de $\brac_+$ dans $\dend$ est un idéal.

\vspace{0.3cm}

On rappelle qu'un $\smod$-module $E$ est une collection de
$\sym_n$-modules $E(n)$ pour $n\geq 1$. Si $e_m\in\PP(m)$,
$f_n\in\PP(n)$ et $i\in\{1,\dots,m\}$, on note $e_m\circ_i f_n$ la
composition de $f_n$ à la place $i$ de $e_m$. On renvoie le lecteur
peu familier avec les opérades à la référence standard \cite{may97}.

\begin{defi}
  Un idéal à gauche dans une opérade $\PP$ est un sous-$\smod$-module
  $\JJ$ tel que pour tous $f_n\in\JJ(n)$, $e_m\in\PP(m)$ et
  $i\in\{1,\dots,m\}$, on ait $e_m\circ_i f_n\in \JJ$.
\end{defi}
On peut rapprocher cette notion de celle d'idéal à gauche dans un
anneau associatif. En particulier, un idéal est aussi un idéal à gauche.

Voici un critère simple pour montrer qu'un idéal à gauche est un
idéal. La démonstration est une conséquence facile des axiomes
d'opérades et est laissée au lecteur.
\begin{lemm}
  \label{critere}
  Soient $\PP$ une opérade et $\JJ$ un idéal à gauche de $\PP$.
  Soient $E$ un système de générateurs de $\PP$ comme
  opérade et $F$ un système de générateurs de $\JJ$
  comme idéal à gauche. Alors $\JJ$ coïncide avec l'idéal
  engendré par $F$ si et seulement si, pour tout $e_n$ dans
  $E(n)$, tout $f_m$ dans $F(m)$ et tout $i\in\{1,\dots,m\}$, on a
  $f_m\circ_i e_n \in\JJ$.
\end{lemm}

\vspace{0.5cm}

Soit $\JJ$ l'idéal à gauche engendré par l'image de $\brac_+$ dans
$\dend$. On aura besoin de deux lemmes préliminaires sur le
$\smod$-module quotient de $\dend$ par $\JJ$.

\begin{lemm}
    \label{droite-gauche}
    Pour tout entier non nul $n$, on a dans $\dend(n)$ modulo $\JJ$,
    \begin{equation}
      \mont(x_n,\dots,x_1)\equiv \desc(x_1,\dots,x_n).
    \end{equation}
  \end{lemm}
  
  \begin{preuve}
    Par récurrence sur $n$. La formule est vraie pour $n=1$. Soit
    maintenant $n\geq 2$ fixé et supposons la formule vérifiée pour
    tous les entiers strictement inférieurs à $n$. On a
    $\mont(x_n,\dots,x_1)\equiv  x_n\prec \, \mont(x_{n-1},\dots,x_1)$, donc
    par hypothèse de récurrence, $\mont(x_n,\dots,x_1)\equiv  x_n\prec\,
    \desc (x_1,\dots,x_{n-1})$. Ceci se réécrit modulo $\JJ$, en
    ajoutant $\psi(\cor_{n-1}(x_n;x_1,\dots,x_{n-1}))$,
    \begin{equation}
      \label{pour-sous-lemme}
      \sum_{i=1}^{n-1} 
      (-1)^i \mont(x_1,\dots,x_{i}) \succ x_n \prec  
      \, \desc(x_{i+1},\dots,x_{n-1}).
    \end{equation}
    On va utiliser le sous-lemme suivant dans une seconde
    récurrence décroissante sur $k$.
    \begin{sublemm}
      Pour $k \geq 2$, modulo $\JJ$, on a
      \begin{multline}
        \sum_{i=1}^{k} 
        (-1)^i \mont(x_1,\dots,x_{i}) \succ z \prec  
        \, \desc(x_{i+1},\dots,x_{k})\equiv \\
        \sum_{i=1}^{k-1} 
        (-1)^i \mont(y_1,\dots,y_{i}) \succ z \prec  
        \, \desc(y_{i+1},\dots,y_{k-1}).
      \end{multline}
      où $y_{1}=x_{1}\succ x_{2}$ et $y_i=x_{i+1}$ pour $2\leq i\leq
      k-1$.
    \end{sublemm}
    \textsc{Preuve du sous-lemme} : 
    La différence entre les deux membres est
    \begin{equation*}
      x_1\succ(\psi(\cor_{k-1}(z;x_2,\dots,x_k))),
    \end{equation*}
    qui appartient bien à $\JJ$. \textsc{C.Q.F.D.}

    \vspace{0.3cm}
    
    Par application répétée du sous-lemme à la formule
    (\ref{pour-sous-lemme}), on la réécrit $ \desc(x_1,\dots,x_n)$, ce
    qui termine la démonstration du lemme.
  \end{preuve}
  
  \vspace{0.3cm} 
  
  On note $\plig(p,q)$ l'ensemble des permutations de
  $\{1,\dots,{p+q}\}$ telles que $\sigma(p)<\dots<\sigma({1})$ et
  $\sigma({p+1})<\dots<\sigma({p+q})$. Ces permutations sont des
  battages où on a renversé l'ordre d'une des deux parties avant de
  battre.
  
  Le lemme suivant permet de réécrire toute opération de $\dend$
  modulo $\JJ$ comme une somme sur certains ensembles de battages.

  \begin{lemm}
    \label{battagedeux}
    Pour tous entiers strictement positifs $p,q$, on a modulo $\JJ$,
    \begin{equation*}
      \desc(x_p,\dots,x_1)\succ z\prec \,\mont(x_{p+q},\dots,x_{p+1})
      \equiv 
      \sum_{\sigma\in\plig(p,q)}\desc(x_{\sigma(1)},\dots,x_{\sigma(p+q)},z).
    \end{equation*}
  \end{lemm}
  \begin{preuve}
    Par récurrence sur le couple $(p,q)$. La formule est vraie pour
    $(p,q)=(1,1)$. Par le lemme \ref{droite-gauche}, on a, modulo $\JJ$,
    \begin{equation*}
      z\prec \,\mont(x_{p+q},\dots,x_{p+1})\equiv 
      \desc(x_{p+1},\dots,x_{p+q})\succ z\equiv 
      \mont(x_{p+q},\dots,x_{p+1})\succ z.
    \end{equation*}
    On en déduit que
    \begin{multline*}
      \desc(x_p,\dots,x_1)\succ z\prec
      \,\mont(x_{p+q},\dots,x_{p+1})\equiv \\
      \big{(}\desc(x_p,\dots,x_{2})\succ x_1 \prec
      \,\mont(x_{p+q},\dots,x_{p+1})\big{)}\succ z\\
      +\big{(}\desc(x_p,\dots,x_1)\succ x_{p+q} \prec
      \,\mont(x_{p+q-1},\dots,x_{p+1})\big{)}\succ z.
    \end{multline*}
    Ceci permet de conclure sans difficulté.
  \end{preuve}
  
  On va maintenant démontrer la proposition suivante :
  
  \begin{prop}
    \label{idealgauche-bras}
    L'idéal à gauche $\JJ$ engendré par l'image de $\brac_+$ dans
    $\dend$ est un idéal.
  \end{prop}
  \begin{preuve}
    D'après le lemme \ref{critere}, et en tenant compte de
    l'automorphisme de $\dend$ donné par la symétrie des arbres plans,
    il suffit de vérifier que $\psi(\cor_n)\circ_j \prec$ est dans
    $\JJ$ pour tout $j$. On distingue deux cas, selon que la
    composition est dans la racine ou dans une feuille de la corolle,
    c'est-à-dire $j=1$ ou $j>1$.
    
    On traite d'abord le cas de la composition à la racine,
    c'est-à-dire
    \begin{equation}
      \psi(\cor_n(z;x_1,\dots,x_n))\circ_z(w\prec t),
    \end{equation}
    soit encore
    \begin{equation}
      \sum_{i=0}^{n} 
      (-1)^i \mont(x_1,\dots,x_i) \succ \big{(}(w\prec t) \prec  
      \, \desc(x_{i+1},\dots,x_{n})\big{)},
    \end{equation}
    ce qui est égal modulo $\JJ$, par application du lemme
    \ref{droite-gauche}, à
    \begin{equation}
      \sum_{i=0}^{n} 
      (-1)^i \desc(x_1,\dots,x_i)\succ\big{(} t\succ w\prec\,
    \mont(x_{i+1},\dots,x_{n})\big{)}.
  \end{equation}
  On utilise deux fois le lemme \ref{battagedeux} pour réécrire ceci
  comme une somme sur $i$ et sur un ensemble de battages. Chaque
  battage s'obtient alors pour deux indices $i$ successifs, donc
  contribue par zéro à la somme, qui est donc nulle.
  
  On traite de façon similaire le cas de la composition dans une
  feuille, en se ramenant à une somme de battage à l'aide des lemmes
  \ref{droite-gauche} et \ref{battagedeux}.
\end{preuve}

\section{Bigèbres dendriformes}

On se place dans cette section et dans toute la suite sur un corps
$\kk$. On définit des notions d'algèbre dendriforme unitaire et de
bigèbre dendriforme. Après avoir rappelé que, d'après \cite{ronco},
les éléments primitifs d'une bigèbre dendriforme $D$ forment une
sous-algèbre brace de $D_{\brac}$, on énonce quelques lemmes sur les
bigèbres dendriformes qui seront utiles dans la suite. On donne une
première famille d'exemples de bigèbres dendriformes : les algèbres
dendriformes libres.

\subsection{Définitions}

On introduit d'abord une notion d'algèbre dendriforme unitaire, voir
\cite[p. 32, (5.5.4)]{dialgebra}. Je remercie Marìa Ronco pour m'avoir
signalé une erreur importante dans une version précédente de cette
définition.

\begin{defi}
  On appelle algèbre dendriforme unitaire un espace vectoriel $W$ muni
  d'une décomposition $W=\kk\cdot 1\oplus V$ et de deux applications
  $\prec,\,\succ : V\otimes W+W\otimes V\to V$ telles que
  $(V,\prec,\succ)$ soit une algèbre dendriforme et que, pour tout
  $x\in V$, on ait
  \begin{equation}
    \label{unite}
    \begin{split}
      1\prec x=x\succ 1=0, \\
      x\prec 1=1\succ x=x.
    \end{split}
  \end{equation}
\end{defi}

Si $D$ est une algèbre dendriforme, alors $\widetilde{D}:=\kk\cdot 1\oplus
D$ est naturellement munie d'une structure d'algèbre dendriforme
unitaire. On associe une algèbre associative unitaire augmentée à une
algèbre dendriforme unitaire en posant $x*y=x\prec y+x\succ y$ et $1*1=1$.

A l'aide de la remarque qui identifie les $\zin$-algèbres aux
algèbres dendriformes symétriques, on définit de
manière analogue la notion de $\zin$-algèbre unitaire.

On appelle algèbre dendriforme unitaire filtrée (resp.
graduée) une algèbre dendriforme unitaire $W$ munie d'une
filtration (resp. d'une graduation) $(W_n)_{n\geq 0}$ telle que $1\in
W_0$ et si $x\in W_p$ et $y\in W_q$, alors $x\prec y\in W_{p+q}$ et
$x\succ y \in W_{p+q}$. L'espace gradué associé à une
algèbre dendriforme unitaire filtrée est naturellement une
algèbre dendriforme unitaire graduée.

Si $W$ est une algèbre dendriforme unitaire, on note $W^+$ le second
facteur de la décomposition $W=\kk\cdot 1\oplus W^+$. Un morphisme
d'algèbres dendriformes unitaires $f :D_1 \to D_2$ est une application qui
respecte les décompositions et telle que $f(1)=1$ et $f:D^+_1\to
D^+_2$ soit un morphisme dendriforme. Ceci définit la catégorie des
algèbres dendriformes unitaires, notée $\dend-\alg-u$. Le lemme
suivant est évident.

\begin{lemm}
  \label{dend-unit}
  Le foncteur $D \mapsto \widetilde{D}$ de $\dend-\alg$ dans
  $\dend-\alg-u$ est une équivalence de catégories, de quasi-inverse
  $W\mapsto W^+$.
\end{lemm}

\vspace{0.5cm}

On introduit maintenant une notion de bigèbre dendriforme. Cette
défi-nition est voisine de celle donnée par Ronco dans
\cite{ronco} sous le nom d'algèbre de Hopf dendriforme. La
différence réside dans la présence d'une unité, qui
simplifie quelque peu la forme des axiomes et de certaines
démons-trations, au prix d'un certain abus de notation.

\begin{defi}
  \label{defi-bigeb}
  On appelle bigèbre dendriforme une algèbre dendriforme unitaire
  $W=\kk\cdot 1\oplus V$ munie d'un coproduit coassociatif $\Delta
  :W\rightarrow W\otimes W$ vérifiant les conditions suivantes :
  \begin{itemize}
  \item $\Delta(1)=1\otimes 1$,
  \item la projection sur $\kk\cdot 1$ parallèlement à $V$ est une
    coünité.
  \end{itemize}
  Ces conditions entraînent, pour $v\in V$, $\Delta(v) \in 1\otimes
  v+v\otimes 1+V\otimes V$.
  \begin{itemize}
  \item Le coproduit est compatible aux produits dendriformes, au sens
    suivant : pour tous $x,y \in V$, avec la notation de Sweedler,
  \begin{equation}
    \label{gauche-coproduit}
    \Delta(x\succ y)=(x_{(1)}*y_{(1)})\otimes(x_{(2)}\succ y_{(2)})
     -(x * y)\otimes (1\succ 1) +(x\succ y)\otimes 1,
  \end{equation}
  et
  \begin{equation}
    \label{droite-coproduit}
    \Delta(x\prec y)=(x_{(1)}*y_{(1)})\otimes(x_{(2)}\prec y_{(2)})
    -(x * y)\otimes (1\prec 1) +(x\prec y)\otimes 1,
  \end{equation}
  \end{itemize}
  où la soustraction des termes en $1\succ 1$ et $1\prec 1$ est une
  notation commode pour signifier que dans la sommation de Sweedler,
  on doit omettre l'unique terme correspondant à $x_{(2)}=y_{(2)}=1$.
\end{defi}

On utilisera par la suite sans davantage de commentaires cette
convention sur les termes non-existants en $1\succ 1$ et $1\prec 1$.

\vspace{0.3cm}

Si $(W,\prec,\succ,\Delta)$ est une bigèbre dendriforme, alors $W$
est en particulier une bigèbre pour le produit associatif $*$ et
le coproduit $\Delta$. 

\vspace{0.3cm}

On appelle bigèbre dendriforme filtrée (resp. graduée) une bigèbre
dendriforme munie d'une filtration (resp. d'une graduation) qui en
fait à la fois une algèbre dendriforme unitaire filtrée (resp.
graduée) et une cogèbre filtrée (resp. graduée). L'espace gradué
associé à une bigèbre dendriforme filtrée est naturellement une
bigèbre dendriforme graduée.

\subsection{Éléments primitifs d'une bigèbre dendriforme}
Si $C$ est une cogèbre, on note $\prim(C)$ le sous-espace des éléments
primitifs de $C$. On rappelle le fait remarquable, établi par Ronco
\cite{ronco}, que les éléments primitifs d'une bigèbre dendriforme $D$
forment une sous-algèbre brace de $D^+_{\brac}$.

\vspace {0.3cm}

On aura besoin du lemme suivant dans les sections 5 et 6.
\begin{lemm}
  \label{coprod-mont}
  Soient $(W,\prec,\succ,\Delta)$ une bigèbre dendriforme et
  $x_1,\dots,x_n$ des éléments primitifs de $W$. Alors
  \begin{equation}
    \Delta(\mont(x_1,\dots,x_n))=\sum_{i=0}^n \mont(x_1,\dots,x_i)
    \otimes\mont(x_{i+1},\dots,x_n),
  \end{equation}
  où on convient que $\mont()=1$.
\end{lemm}
\begin{preuve}
  Par récurrence, à l'aide des formules (\ref{gauche-coproduit}) et
  (\ref{droite-coproduit}), voir \cite[2.7]{ronco}, pour plus de
  détails.
\end{preuve}

\vspace {0.3cm}

\begin{prop}
  \label{primitif}
  Si $W=\kk\cdot 1\oplus V$ est une bigèbre dendriforme, l'ensemble
  $\prim(W)$ des éléments primitifs de $W$ est une sous-algèbre
  brace de $V_{\brac}$.
\end{prop}
\begin{preuve}
  Il résulte aussitôt de la définition que $\prim W$ est
  contenu dans $V$. On renvoie à \cite[Prop 2.8]{ronco} pour la
  démonstration du fait que $\prim W$ est une sous-algèbre brace.
\end{preuve}

\subsection{Lemmes}

Les lemmes suivants, qui sont des conséquences directes des axiomes de
bigèbre dendriforme, jouent un rôle crucial pour la suite.
\begin{lemm}
  \label{coideal}
  Si $W=\kk\cdot 1\oplus V$ est une bigèbre dendriforme et $I\subset
  V$ un coidéal de $W$, alors l'idéal dendriforme de $V$ engendré par
  $I$ est un coidéal de $W$.
\end{lemm}
\begin{preuve}
  Soit $I$ un coidéal de $W$ contenu dans $V$. Il résulte des formules
  (\ref{gauche-coproduit}) et (\ref{droite-coproduit}) que
  $\mathfrak{s}(I)=I+V\prec I+I\prec V+V\succ I+I\succ V$ est aussi un
  coidéal de $W$ contenu dans $V$. On définit par récurrence
  \begin{equation*}
    \mathfrak{s}^n(I):=\mathfrak{s}(\mathfrak{s}^{n-1}(I)).
  \end{equation*}
  Les coidéaux $\mathfrak{s}^n(I)$ sont inclus dans $V$ et forment une
  suite croissante. Soit alors $\mathfrak{S}(I)$ l'union des
  $\mathfrak{s}^n(I)$ pour $n\geq 1$. C'est encore un coidéal de $W$
  et il est clair que c'est l'idéal dendriforme de $V$ engendré par
  $I$. Ceci termine la démonstration.
\end{preuve}

\vspace{0.5cm}

\begin{lemm}
  \label{critere-coproduit}
  Soient $D_1=\kk\cdot 1\oplus D^+_1,D_2=\kk\cdot 1\oplus D^+_2$ deux
  bigèbres dendriformes et $f : D_1 \to D_2$ un morphisme d'algèbres
  dendriformes unitaires. Si $D_1^+$ est engendrée par $V$ comme algèbre
  dendriforme et si $\Delta\circ f=(f\otimes f)\circ\Delta$
  sur $V$, alors $f$ est un morphisme de bigèbres dendriformes.
\end{lemm}
\begin{preuve}
  On a clairement $\Delta(f(1))=(f\otimes f) \Delta(1)$. Considérons
  le sous-espace
  \begin{equation*}
    K:=\{x\in D_1^+ \mid  \Delta(f(x))=(f\otimes f) \Delta(x)\}.
  \end{equation*}
  Par hypothèse, $K$ contient $V$. On montre
  sans difficulté à l'aide des formules
  (\ref{gauche-coproduit}) et (\ref{droite-coproduit}) que $K$ est une
  sous-algèbre dendriforme de $D_1^+$. Comme $D_1^+$ est engendrée
  par $V$, on a $K=D_1^+$, ce qui termine la démonstration.
\end{preuve}

\subsection{Algèbres dendriformes libres}
On montre que les algèbres dendriformes libres sont des bigèbres
dendriformes. Ce résultat a également été obtenu par Marìa Ronco
\cite{ronco}. On en présente ici une démonstration légèrement
différente qui utilise la présence d'une unité. On considère l'algèbre
dendriforme unitaire $\widetilde{D}(V):=\kk\cdot 1\oplus D(V)$ où
$D(V)$ est l'algèbre dendriforme libre sur $V$.

Le coproduit est celui introduit par Loday et Ronco dans
\cite{loday-ronco} qui, \textit{stricto sensu}, ne traite que le cas
de l'algèbre libre sur un générateur. L'extension des
résultats de \cite{loday-ronco} au cas général est
toutefois immédiate, voir \cite{ronco} pour plus de détails.
On obtient ainsi la description suivante du coproduit $\Delta$.

Pour $x,y\in \widetilde{D}(V)$ et $v\in V$, on pose d'abord
\begin{equation*}
    x \viv y=x\succ v \prec y.
\end{equation*}
Alors, le coproduit $\Delta$ vérifie, pour $x,y\in
\widetilde{D}(V)$ et $v\in V$,
\begin{equation}
\label{coprod-vee}
    \Delta(x\viv y)=(x\viv y)\otimes 1+(x_{(1)}* y_{(1)})
    \otimes(x_{(2)}\viv y_{(2)}),
\end{equation}
et il est entièrement déterminé par ces relations et la
condition $\Delta(1)=1\otimes 1$. On remarque que le second terme du
membre de droite contient $1\otimes (x\viv y)$.

\begin{prop}
  \label{dendri-hopf}
  L'algèbre dendriforme unitaire $\widetilde{D}(V)$, munie du
  coproduit $\Delta$, est une bigèbre dendriforme.
\end{prop}
\begin{preuve}
  On déduit de l'identité $ x\succ(y\succ z)=(x*y)\succ z$
  dans $D(V)$ la formule suivante, pour $x,w,z\in D(V)$ et $v\in V$ :
  \begin{equation}
    x\succ(w\viv z)=(x*w)\viv z.
  \end{equation}
  Elle est encore valable pour $x,w,z\in\widetilde{D}(V)$.

  Soient maintenant $x,y\in D(V)$. Comme on suppose que $y\in D(V)$, on peut
  écrire de manière unique $y=w\viv z$ avec $w,z\in
  \widetilde{D}(V)$ et $v\in V$. On a alors
  \begin{multline*}
    \Delta(x\succ y)=\Delta(x\succ(w\viv z))=\Delta((x*w)\viv z)\\
    =((x*w)\viv z)\otimes 1+((x*w)_{(1)}*z_{(1)})
    \otimes((x*w)_{(2)}\viv z_{(2)}).
  \end{multline*}
  Comme $\Delta$ est un morphisme d'algèbres associatives, ceci est
   égal à
  \begin{equation}
  ((x*w)\viv z)\otimes 1+(x_{(1)}*w_{(1)}*z_{(1)})
    \otimes((x_{(2)}*w_{(2)})\viv z_{(2)}).
  \end{equation}
  D'autre part, en utilisant la convention des termes fantômes en
  $1\succ 1$ (voir Def. \ref{defi-bigeb}), on a 
  \begin{multline}
    (x_{(1)}*y_{(1)})\otimes(x_{(2)}\succ y_{(2)})
    -(x * y)\otimes (1\succ 1) +(x\succ y)\otimes 1 =\\
    (x_{(1)}*(w\viv z)_{(1)})\otimes(x_{(2)}\succ (w\viv
    z)_{(2)})-(x * (w\viv
    z))\otimes (1\succ 1)+(x\succ (w\viv z))\otimes 1 \\
    =( x*(w\viv z))\otimes (1 \succ 1)+(x_{(1)}*w_{(1)}*z_{(1)})
    \otimes(x_{(2)}\succ(w_{(2)}\viv z_{(2)}))\\
    -(x * (w\viv
    z))\otimes (1\succ 1)+(x\succ (w\viv z))\otimes 1 \\
    =((x*w)\viv z)\otimes 1+(x_{(1)}*w_{(1)}*z_{(1)})
    \otimes((x_{(2)}*w_{(2)})\viv z_{(2)}).
  \end{multline}
  
  Ceci termine la preuve de la formule (\ref{gauche-coproduit}). La
  formule (\ref{droite-coproduit}) s'en déduit par l'automorphisme
  de symétrie gauche-droite des arbres plans.
\end{preuve}
\section{Foncteurs adjoints et suites exactes courtes}

\subsection{Foncteurs adjoints}

La définition de l'algèbre dendriforme enveloppante passe d'abord par
celle de l'adjoint à gauche du foncteur d'oubli $()_{\brac}$ qui
associe à une algèbre dendriforme l'algèbre brace sous-jacente. La
proposition suivante a pour vocation essentielle de rappeler la
construction générale de ce type de foncteur adjoint, qui est sans
difficulté et sans doute bien connue.
\begin{prop}
  Soit $\gamma : \PP\to\QQ$ un morphisme d'opérades. Le foncteur
  d'oubli $\Gamma : \QQ-\alg\to \PP-\alg$ admet un adjoint à gauche,
  noté $U_{\Gamma}$.
\end{prop}
\begin{preuve}
  Soit $P$ une $\PP$-algèbre, on définit $U_{\Gamma}(P)$ comme
  le quotient de la $\QQ$-algèbre libre sur $P$, notée
  $F_{\QQ}P$, par le $\QQ$-idéal engendré par les relations
  \begin{equation}
    \gamma(e_n)(p_1,\dots,p_n)-e_n(p_1,\dots,p_n)
  \end{equation}
  pour tous $n\geq 1$, $e_n\in\PP(n)$ et $p_1,\dots,p_n$ dans $P$. On
  observe qu'il revient au même de quotienter par ces relations pour
  des $e_n$ décrivant seulement un ensemble $E$ de générateurs de $\PP$.
  
  Par construction, l'application naturelle de $P$ dans $F_{\QQ}P$
  fournit un $\PP$-morphisme $\tau_{P}$ de $P$ dans $U_{\Gamma}(P)$.
  Il faut vérifier qu'il a bien la propriété universelle voulue.
  
  Soient $D$ une $\QQ$-algèbre et $\phi : P\rightarrow \Gamma(D)$
  un $\PP$-morphisme. Alors $\phi$ se prolonge en un $\QQ$-morphisme
  $\Phi : F_{\QQ}P \rightarrow D$. Comme $\phi$ est un un
  $\PP$-morphisme, $\Phi$ se factorise en un $\QQ$-morphisme
  $\overline{\Phi} : U_{\Gamma}(P)\rightarrow D$. On a clairement
  $\overline{\Phi}\circ\tau_{P}=\phi$. Un $\QQ$-morphisme qui
  vérifie cette relation est unique, car $U_{\Gamma}(P)$ est
  engendrée comme $\QQ$-algèbre par $\tau_{P}(P)$.

  On a donc une bijection, donnée par la composition avec $\tau_P$,
  \begin{equation*}
    \hom_{\PP-\alg}(P,\Gamma(D))\simeq\hom_{\QQ-\alg}(U_{\Gamma}(P),D). 
  \end{equation*}

  On vérifie sans peine que l'on obtient ainsi une équivalence de
  bifoncteurs. Ceci démontre la proposition.
\end{preuve}

\vspace{0.5cm}

La cas particulier qui nous intéresse ici est le suivant :

\begin{coro}
  Le foncteur $()_{\brac}$ de $\dend-\alg$ dans $\brac-\alg$ admet
  un adjoint à gauche, noté $\UU$.
\end{coro}

Explicitement, si $B$ est une algèbre brace, $\UU(B)$ est le
quotient de l'algèbre dendriforme libre $D(B)$ par l'idéal
dendriforme engendré par les éléments
\begin{equation}
  \label{rela-quot}
  \psi(\cor_n)(x_1,x_2,\dots,x_{n+1})-\cor_n(x_1,x_2,\dots,x_{n+1}),
\end{equation}
pour $n\geq 2$ et $x_1,\dots,x_{n+1}\in B$.
  
\vspace{0.2cm}

Le cas des algèbres braces libres $B(V)$ est
particulièrement simple. On le développe ici pour un usage
ultérieur.
\begin{lemm}
  \label{explicite-libre}
  Les algèbres dendriformes $\UU(B(V))$ et $D(V)$ sont isomorphes.
\end{lemm}
\begin{preuve}
  Le foncteur d'oubli $\dend-\alg\rightarrow \vect$ est la
  composée des foncteurs d'oubli $\dend-\alg\rightarrow
  \brac-\alg$ et $\brac-\alg\rightarrow \vect$. Comme l'adjoint
  d'un foncteur composé est la composée des foncteurs
  adjoints, et que pour toute opérade $\PP$, le foncteur
  $\PP$-algèbre libre est adjoint au foncteur d'oubli de
  $\PP-\alg\rightarrow \vect$, on obtient le résultat.
\end{preuve}

\subsection{Suites exactes courtes d'opérades}
Soit $P$ une $\PP$-algèbre. Alors, l'ensemble $\JJ$ des opérations de
$\PP$ qui sont nulles dans $P$, c'est-à-dire le noyau du morphisme
$\PP\to\aut(P)$, forme un idéal (bilatère) au sens des opérades. On
cherche ici sous quelles conditions la nullité de certaines opérations
de $\PP$ sur $P$ peut se déduire de leur nullité sur un sous-espace
$V$ de $P$ engendrant $P$ comme $\PP$-algèbre.

\vspace{0.3cm}

On omettra ici les démonstrations des deux prochains lemmes, qui sont
des exercices de manipulation des axiomes d'opérades, sans surprises.
Le lemme suivant justifie l'introduction de la notion d'idéal à gauche.
\begin{lemm}
\label{id-gauche}
  Soit $P$ une $\PP$-algèbre et $V$ un sous-espace de $P$.
  L'ensemble des opérations de $\PP$ nulles sur $V$ forme un
  idéal à gauche de $\PP$.
\end{lemm}

\begin{lemm}
\label{ideal-nul}
  Soit $P$ une $\PP$-algèbre, $V$ un sous-espace de $P$ qui
  engendre $P$ comme $\PP$-algèbre et $\JJ$ un idéal de $\PP$.
  Si les opérations de $\JJ$ sont nulles sur $V$, elles sont
  nulles sur $P$.
\end{lemm}

\vspace{0.3cm}

Soient $\PP$ une opérade augmentée, $\PP_+$ son idéal d'augmentation
et $\gamma : \PP\to\QQ$ un morphisme d'opérades. On note comme
ci-dessus $U_{\Gamma}$ l'adjoint du foncteur $\Gamma$ de $\QQ-\alg$
dans $\PP-\alg$ induit par $\gamma$. On a alors la

\begin{prop}
  \label{idealgauche-ideal}
  On suppose que $\JJ$, l'idéal à gauche de $\QQ$ engendré par
  $\gamma(\PP_+)$ est un idéal. On note $\CC$ l'opérade quotient de
  $\QQ$ par $\JJ$. Si $P$ est une $\PP$-algèbre annulée par $\PP_+$,
  alors $U_{\Gamma}(P)$ est une $\CC$-algèbre, isomorphe à la
  $\CC$-algèbre libre sur $P$.
\end{prop}
\begin{preuve}
  Soit $\tau$ l'inclusion de $P$ dans $U_{\Gamma}(P)$. Il résulte de
  l'hypothèse et de la définition de $U_{\Gamma}(P)$ que les
  opérations $\gamma(\PP_+)$ sont nulles sur $\tau(P)$. Donc, d'après
  le lemme \ref{id-gauche}, les opérations de $\JJ$ sont nulles sur
  $\tau(P)$. Comme $U_{\Gamma}(P)$ est engendrée comme $\QQ$-algèbre
  par $\tau(P)$ et comme, par hypothèse, $\JJ$ est un idéal, il
  résulte du lemme \ref{ideal-nul} que les opérations de $\JJ$ sont
  nulles sur $U_{\Gamma}(P)$. Par conséquent, $U_{\Gamma}(P)$ est une
  $\CC$-algèbre.
  
  Il reste à montrer que la $\CC$-algèbre $U_{\Gamma}(P)$, possède la
  propriété universelle de la $\CC$-algèbre libre sur $P$. Soit $Z$
  une $\CC$-algèbre et $\phi:P\to Z$ un morphisme d'espaces
  vectoriels. Comme $Z$ est aussi une $\QQ$-algèbre, il existe un
  unique $\QQ$-morphisme $\widetilde{\phi}:F_{\QQ}P\to Z$ prolongeant
  $\phi$. Comme $Z$ est une $\CC$-algèbre, ce morphisme passe au
  quotient pour donner un $\QQ$-morphisme
  $\overline{\phi}:U_{\Gamma}(P)\to Z$ tel que $\overline{\phi} \circ
  \tau=\phi$. Comme $U_{\Gamma}(P)$ et $Z$ sont des $\CC$-algèbres,
  $\overline{\phi}$ est un $\CC$-morphisme. Enfin, un tel morphisme
  est unique, car $P$ engendre $U_{\Gamma}(P)$ comme $\QQ$-algèbre,
  donc aussi comme $\CC$-algèbre.
\end{preuve}

\section{Algèbre dendriforme enveloppante}
\subsection{Définition et coproduit}
On définit l'\textit{algèbre dendriforme enveloppante} d'une algèbre
brace $B$, notée $\widetilde{\UU}(B)$, comme l'algèbre dendriforme
unitaire associée à $\UU(B)$.

\begin{prop}
  \label{coproduit-quotient}
  $\widetilde{\UU}(B)$ est une bigèbre dendriforme.
\end{prop}
\begin{preuve}
  Comme $\prim(\widetilde{D}(B))$ est une sous $\brac$-algèbre de
  $D(B)$, les éléments apparaissant dans (\ref{rela-quot})
  appartiennent à $\prim(\widetilde{D}(B))$, et donc forment un
  coidéal de $\widetilde{D}(B)$. Par le lemme \ref{coideal}, l'idéal
  dendriforme de $D(B)$ engendré par ces éléments est encore un coidéal de
  $\widetilde{D}(B)$, donc la structure de bigèbre dendriforme de
  $\widetilde{D}(B)$ passe au quotient.
\end{preuve}

\vspace{0.5cm} 

On note $B(V)$ l'algèbre brace libre sur un espace vectoriel $V$.

Par le lemme \ref{explicite-libre}, l'algèbre dendriforme enveloppante
$\widetilde{\UU}(B(V))$ est isomorphe à $\widetilde{D}(V)$. On
retrouve ainsi la bigèbre dendriforme des arbres binaires plans, dans
la terminologie de \cite{loday-ronco}, voir le paragraphe 3.4.

\subsection{Cas des algèbres braces de produits nuls}

On appelle algèbre brace triviale une algèbre brace dont tous les
produits $\{\,?\,\mid\, ?\,\}_n$ sont nuls.

On donne une description en termes d'objets classiques de
l'algèbre dendriforme enveloppante d'une algèbre brace
triviale.

\vspace{0.5cm} 

On note $T^{c}(V)$ la cogèbre tensorielle sur $V$, voir
\cite[1.4,1.5]{reutenauer} pour la définition. On rappelle que
$T^{c}(V)$ possède une structure de bigèbre, où le
coproduit est donné par la déconcatenation et le produit
commutatif par les battages, et telle que $\prim(T^{c}(V))=V$, cf
\cite[loc. cit.]{reutenauer}. On rappelle deux propriétés
classiques : $T^{c}(V)$ est la cogèbre coassociative connexe colibre sur
$V$ et $T^{c}(V)$ est duale au sens gradué de l'algèbre
tensorielle $T(V)$.

On note $\widetilde{Z}(V)$ la $\zin$-algèbre unitaire $\kk\cdot
1\oplus Z(V)$. D'après \cite[7.1]{dialgebra},
$\widetilde{Z}(V)_{\com}$ s'identifie comme algèbre commutative à
$T^{c}(V)$. Par conséquent, $T^{c}(V)$ est naturellement une algèbre
dendriforme symétrique et $\widetilde{Z}(V)_{\com}$ une bigèbre.

Soit maintenant $V$ une algèbre brace triviale. On a alors la
description suivante de $\widetilde{\UU}(V)$ :

\begin{prop}
  Soit $V$ une algèbre brace triviale.
\begin{itemize}
\item Il existe un unique isomorphisme d'algèbres dendriformes
  unitaires $\theta : \widetilde{\UU}(V) \to \widetilde{Z}(V)$ qui est
  l'identité sur $V$.
\item De plus, $\theta$ est un isomorphisme de cogèbres.
\end{itemize}
\end{prop}
\begin{preuve}
  Par la proposition \ref{idealgauche-bras}, on peut appliquer la
  proposition \ref{idealgauche-ideal} au morphisme $\psi :
  \brac\to\dend$. Par la proposition \ref{quotient-zin}, l'opérade
  quotient est $\zin$. Par la proposition \ref{idealgauche-ideal}, il
  existe donc un unique isomorphisme dendriforme de $\UU(V)$ dans
  $Z(V)$ qui est l'identité sur $V$. Par le lemme \ref{dend-unit}, ce
  morphisme se relève de manière unique en un isomorphisme d'algèbres
  dendriformes unitaires de $\widetilde{\UU}(V)$ dans
  $\widetilde{Z}(V)$. On obtient la première assertion.
  
  Il reste à vérifier que le coproduit sur $\widetilde{\UU}(V)$ défini
  par la proposition \ref{coproduit-quotient} coïncide avec le
  coproduit de déconcaténation sur $\widetilde{Z}(V)=T^{c}(V)$.
  Explicitons pour cela l'unique $\zin$-morphisme de $Z(V)$ dans
  $\UU(V)$ qui fixe $V$. On montre sans difficulté que l'image du
  tenseur $x_1\otimes x_2\otimes\dots\otimes x_n$ est la classe de
  l'élément $\mont(x_1,\dots,x_n)$.
  
  Il résulte alors du lemme \ref{coprod-mont} que le coproduit sur
  $\widetilde{\UU}(V)$ correspond par l'isomorphisme à la
  déconcatenation dans $\widetilde{Z}(V)$.
\end{preuve}

\vspace{0.3cm}

On en déduit immédiatement le corollaire suivant, qui a aussi été noté
par Ronco \cite{ronco}.
\begin{coro}
  Il existe sur $T^{c}(V)$ une structure de bigèbre dendriforme dont
  la structure de bigèbre sous-jacente est la structure usuelle de la
  bigèbre des battages.
\end{coro}

\section{Une équivalence de catégories}
On montre dans cette section que le foncteur $\widetilde{\UU}$ est une
équivalence de la catégorie des algèbres braces dans celle des
bigèbres dendriformes connexes, de quasi-inverse le foncteur $\prim$.
Ceci forme un équivalent dendriforme/brace du théorème de
Cartier-Milnor-Moore-Quillen sur les bigèbres cocommutatives connexes
et les algèbres de Lie \cite{milnor-moore,quillen}.

\vspace{0.3cm}

On dit qu'une bigèbre dendriforme sur $\kk$ est connexe si la cogèbre
sous-jacente est connexe au sens de Quillen \cite{quillen},
c'est-à-dire si le coradical est $\kk$ et si la filtration par le
coradical est exhaustive.

Soit $D$ une bigèbre dendriforme et notons $B$ le sous-espace des
éléments primitifs de $D$. D'après la proposition \ref{primitif}, $B$
est une sous-algèbre brace de $D^+_{\brac}$. Par conséquent, d'après
la propriété d'adjonction du foncteur $\UU$, on obtient un morphisme
d'algèbres dendriformes de $\UU(B)$ dans $D^+$, qu'on relève en un
morphisme d'algèbres dendriformes unitaires $\theta$ de
$\widetilde{\UU}(B)$ dans $D$. On a la proposition suivante.

\begin{prop}
  $\theta$ est un morphisme de cogèbres.
\end{prop}
\begin{preuve}
  La restriction de $\theta$ à $B$ vérifie $\Delta\circ\theta=1
  \otimes \theta+\theta\otimes 1=(\theta\otimes\theta)\circ\Delta$.
  Comme $B$ engendre $\UU(B)$ comme algèbre dendriforme, le lemme
  \ref{critere-coproduit} permet de conclure.
\end{preuve}

\vspace{0.3cm}

On rappelle le lemme classique suivant.
\begin{lemm}
\label{injectif-prim}
  Soit $f:C\to C'$ un morphisme de cogèbres. Si $C$ est connexe et si $f$
  restreint aux primitifs de $C$ est injectif, alors $f$ est injectif.
\end{lemm}
\begin{preuve}
  Voir \cite[Appendice B, Prop. 3.2]{quillen} ou \cite[Th. 2.4.11]{abe}.
\end{preuve}

\vspace{0.3cm}

La proposition suivante joue un rôle crucial dans la suite.

\begin{prop}
  \label{cogebre-colibre}
  Soit $B$ une algèbre brace, alors $\widetilde{\UU}(B)$ est isomorphe
  comme cogèbre à $T^c(B)$. Par conséquent, $\prim \widetilde{\UU}(B)=B$.
\end{prop}
\begin{preuve}
  On définit une application linéaire $f$ de $T^c(B)$ dans
  $\widetilde{\UU}(B)$ par $f(x_1\otimes\dots\otimes
  x_n)=\mont(x_1,\dots, x_n)$. Par le lemme \ref{coprod-mont}, $f$ est
  un morphisme de cogèbres. Comme la cogèbre $T^c(B)$ est connexe et
  $f$ injectif sur $B$, le lemme \ref{injectif-prim} entraîne que $f$
  est injectif.
  
  D'autre part, la bigèbre dendriforme $\widetilde{\UU}(B)$ hérite
  comme quotient de $\widetilde{D}(B)$ d'une structure de bigèbre
  dendriforme filtrée. En notant $\JJ$ le noyau du morphisme
  $\widetilde{D}(B)\to\widetilde{\UU}(B)$, c'est-à-dire l'idéal
  dendriforme de $D(B)$ engendré par les $\psi(\cor_n)-\cor_n$, on
  voit aussitôt que l'idéal $gr\JJ$ de $D(B)=gr D(B)$ contient les
  éléments $\psi(\cor_n)$. Par conséquent, $gr\,\widetilde{\UU}(B)$
  est un quotient de $\widetilde{\UU}(B_0)$, où $B_0$ désigne $B$
  munie de la structure brace nulle. Comme $\widetilde{\UU}(B_0)\simeq
  T^c(B)$ comme bigèbre dendriforme, il en résulte que
  $gr\,\widetilde{\UU}(B)$ est engendré comme espace vectoriel par les
  éléments $\mont(x_1,\dots ,x_n)$. Par conséquent,
  $\widetilde{\UU}(B)$ est aussi engendré comme espace vectoriel par
  ces éléments, donc $f$ est surjective.
\end{preuve}

\vspace{0.3cm}

On peut maintenant démontrer la proposition suivante. La preuve
s'inspire de celle donnée par Quillen dans \cite[Appendice B, Th.
4.5]{quillen}.

\begin{prop}
\label{prim-connexe}
  Soit $D$ une bigèbre dendriforme supposée connexe comme cogèbre et
  $B:=\prim(D)$ l'algèbre brace des primitifs de $D$. Alors on a un
  isomorphisme naturel de bigèbres dendriformes $\widetilde{\UU}(B)\to D$.
\end{prop}
\begin{preuve}
  Il s'agit de montrer que le morphisme de bigèbres dendriformes
  naturel $\theta$ de $\widetilde{\UU}(B)$ dans $D$ est bijectif.
  
  Comme $\widetilde{\UU}(B)\simeq T^c(B)$ comme cogèbre,
  $\widetilde{\UU}(B)$ est connexe. Comme $\theta$ est injectif sur
  $B$, $\theta$ est injectif par le lemme \ref{injectif-prim}. Par
  conséquent, en choisissant arbitrairement un supplémentaire de
  $\widetilde{\UU}(B)$ dans $D$, on voit qu'il existe une application
  linéaire $\pi : D\to B$ telle que $\pi\circ\theta $ soit la
  projection canonique $p:\widetilde{\UU}(B)\to B$. Mais, comme
  $\widetilde{\UU}(B)\simeq T^c(B)$ comme cogèbre, elle a la propriété
  universelle de la cogèbre connexe colibre, c'est-à-dire il existe un
  unique morphisme de cogèbres $ g :D\to \widetilde{\UU}(B)$ tel que
  $p\circ g=\pi$. Comme $D$ est connexe par hypothèse et que $g$
  restreint aux primitifs $B$ de $D$ est injectif, $g$ est injectif
  par le lemme \ref{injectif-prim}.
  
  D'autre part, $g\circ \phi$ est un endomorphisme de la cogèbre
  $\widetilde{\UU}(B)$ tel que $g\circ \phi\circ p=p$, donc par la
  propriété universelle, $g\circ \phi=\id$, donc $g$ est aussi
  surjectif. On conclut que $g$ est bijectif d'inverse $\phi$, donc
  $\phi$ est un isomorphisme de bigèbres dendriformes, ce qui termine
  la démonstration.
\end{preuve}

\vspace{0.5cm}

En combinant les propositions \ref{cogebre-colibre} et
\ref{prim-connexe}, on obtient le théorème principal.
\begin{theo}
  Le foncteur $B \mapsto \widetilde{\UU}(B)$ est une équivalence
  entre la catégorie des algèbres braces et celle des
  bigèbres dendriformes connexes, le foncteur quasi-inverse
  étant $ D \mapsto \prim(D)$.
\end{theo}

\vspace{0.5cm}

\textbf{Remerciements.} Je voudrais remercier ici Jim Stasheff pour
son interêt et ses remarques, et mon directeur de thèse
Patrick Polo pour son soutien constant et ses suggestions sur la
dernière partie.

 \bibliographystyle{plain}


Frédéric CHAPOTON\\
Équipe Analyse algébrique, case 82\\
Institut de Mathématiques, Université Pierre et Marie Curie\\
4, Place Jussieu, 75252 Paris Cedex 05, FRANCE\\
Mel. : chapoton$@$math.jussieu.fr

\end{document}